\documentclass[12pt]{article}
\usepackage[cp1251]{inputenc}
\usepackage[russian]{babel}
\usepackage{amsfonts, amsmath, amssymb}
\usepackage{pgf,tikz,url}
\usetikzlibrary{arrows}
\input epsf
\begin{document}

\centerline{\uppercase{\bf Вложения в плоскость графов с вершинами степени 4}}
\bigskip
\centerline{\bf А. Скопенков\footnote{
Московский Физико"=Технический Институт, Независимый Московский Университет. 
Инфо: \url{www.mccme.ru/~skopenko}. 
Благодарю редколлегию сборника <<Математическое Просвещение>>
и особенно Б.Р. Френкина за полезные замечания и обсуждения.
Поддержан РФФИ, грант номер  и~\mbox{15-01-06302}, грантами фонда Саймонса и~стипендией фонда Д. Зимина <<Династия>>.
}}
\bigskip



В этой методической заметке приводится доказательство интересной и важной гипотезы
В.А. Васильева [V] о критерии планарности графа с вершинами степени 4 и некоторой дополнительной
(`крестовой') структурой.
Эта гипотеза относится к числу тех замечательных проблем, возникших в
современной математике, формулировка и доказательство которых доступны школьнику,
знакомому с основами теории графов.
Она была впервые доказана В.О. Мантуровым [M1] 
(это было одним из основных результатов его докторской диссертации [M2]).
Здесь приведено ясное
изложение доказательства,%
\footnote{
Эта заметка написана по материалам цикла задач, представленного на Летней
Конференции Турнира Городов в 2005 г. А. Каибхановым, Д. Пермяковым и автором:
\url{http://olympiads.mccme.ru/lktg/2005/3/index.htm}.
Она выложена в архив в августе 2010.
Поставленный 
вопрос (задача 2) послужил отправной точкой для интересных исследований [F, FM].}
а также пояснения для начинающих, ср. [M3].

\smallskip

{\it Полуребром} графа называется пара из ребра и его концевой вершины
(или, неформально, конец ребра).
{\it Графом с крестовой структурой} (кратко: {\it $X$-графом}) называется
граф, степень каждой вершины которого равна 4, причем для каждой вершины
фиксировано разбиение четырёх входящих в нее полурёбер на две пары.

Вот примеры $X$-графов.
Они изображены на плоскости (с самопересечениями) так, что
при обходе вокруг каждой вершины полурёбра из первой пары и
полурёбра из второй пары чередуются (т.е. проходятся в порядке 1212, а не 1122).

\smallskip
\centerline{
\begin{tabular}{c@{\mbox{\quad}}c@{\mbox{\quad}}c@{\mbox{}}c}
\epsfbox{pict1.19} &
\definecolor{zzttqq}{rgb}{0.6,0.2,0}
\definecolor{qqqqff}{rgb}{0,0,1}
\begin{tikzpicture}[line cap=round,line join=round,>=triangle 45,x=0.5cm,y=0.5cm]
\clip(-2.96,0.98) rectangle (4.16,7.78);
\draw [color=zzttqq] (-1,4)-- (2,4);
\draw (-1,4)-- (0.5,6.6);
\draw (2,4)-- (-1,4);
\draw (2,4)-- (0.48,6.56);
\draw [shift={(-0.96,5.66)}] plot[domain=0.56:4.69,variable=\t]({1*1.7*cos(\t r)+0*1.7*sin(\t r)},{0*1.7*cos(\t r)+1*1.7*sin(\t r)});
\draw [shift={(1.98,5.7)}] plot[domain=-1.56:2.64,variable=\t]({1*1.7*cos(\t r)+0*1.7*sin(\t r)},{0*1.7*cos(\t r)+1*1.7*sin(\t r)});
\draw [shift={(0.5,3.14)}] plot[domain=-3.66:0.52,variable=\t]({1*1.73*cos(\t r)+0*1.73*sin(\t r)},{0*1.73*cos(\t r)+1*1.73*sin(\t r)});
\fill [color=qqqqff] (-1,4) circle (1.5pt);
\fill [color=qqqqff] (2,4) circle (1.5pt);
\fill [color=qqqqff] (0.48,6.56) circle (1.5pt);
\fill [color=qqqqff] (0.48,6.56) circle (1.5pt);
\end{tikzpicture}
&
\epsfbox{pict1.15} &
\definecolor{zzttqq}{rgb}{0.6,0.2,0}
\definecolor{qqqqff}{rgb}{0,0,1}
\begin{tikzpicture}[line cap=round,line join=round,>=triangle 45,x=0.5cm,y=0.5cm]
\clip(-4.3,0.86) rectangle (4.42,8.04);
\draw [color=zzttqq] (-1,4)-- (2,4);
\draw (-1,4)-- (0.5,6.6);
\draw (2,4)-- (-1,4);
\draw (2,4)-- (0.64,6.36);
\draw [shift={(-1,5.64)}] plot[domain=0.72:4.71,variable=\t]({1*1.67*cos(\t r)+0*1.67*sin(\t r)},{0*1.67*cos(\t r)+1*1.67*sin(\t r)});
\draw [shift={(2,5.72)}] plot[domain=-1.57:2.67,variable=\t]({1*1.72*cos(\t r)+0*1.72*sin(\t r)},{0*1.72*cos(\t r)+1*1.72*sin(\t r)});
\draw [shift={(0.5,3.1)}] plot[domain=-3.68:0.54,variable=\t]({1*1.75*cos(\t r)+0*1.75*sin(\t r)},{0*1.75*cos(\t r)+1*1.75*sin(\t r)});
\fill [color=qqqqff] (-1,4) circle (1.5pt);
\fill [color=qqqqff] (2,4) circle (1.5pt);
\end{tikzpicture}\\
Рисунок 1 & Рисунок 2 & Рисунок 3 & Рисунок 4\end{tabular}
}


\smallskip
$X$-граф называется {\it $X$-планарным},
если его можно изобразить без самопересечений на плоскости так,
что при обходе вокруг каждой вершины полурёбра из первой пары и
полурёбра из второй пары чередуются (т.е. проходятся в порядке 1212, а не 1122).

Например, $X$-граф на рис. 1 является $X$-планарным,
ибо он
<<не отличается как $X$-граф>> от изображённого на рис. 2.

Заметим, что не любой планарный $X$-граф $X$-планарен.
Действительно, $X$-графы <<восьмёрка>> (рис. 3) и <<трилистник>> (рис. 4) планарны, но не $X$-планарны.

Общая вершина двух путей на $X$-графе, не имеющих общих рёбер, называется
{\it точкой перекрестья}, если один из этих путей проходит по полурёбрам из одной пары, выходящим из этой вершины.
(Тогда второй путь проходит по полурёбрам из другой пары.)

\smallskip
{\bf Теорема планарности графов c крестовой структурой.}
{\it $X$-граф $X$-планарен тогда и только тогда, когда он не содержит
двух несамопересекающихся циклов без общих рёбер, имеющих ровно одну вершину перекрестья.}

\smallskip
{\it Замечания.}
(a)
Слово <<несамопересекающихся>> можно опустить.
\footnote{
Чтобы обосновать это замечание, покажем, что если в $X$-графе есть
два цикла $K_1$ и $K_2$ ровно с одной точкой перекрестья $A$, то есть и два
несамопересекающихся цикла с тем же свойством.
Для этого рассмотрим
несамопересекающийся цикл $K_1^+$ из рёбер цикла $K_1$, содержащий вершину $A$.
Аналогично определим цикл $K_2^+$. Тогда $K_1^+$ и $K_2^+$ и будут искомыми циклами.
}

(b) Слова <<без общих рёбер, имеющих ровно одну вершину перекрестья>> нельзя заменить на
<<пересекающихся только по одной вершине, являющейся вершиной перекрестья>> (см. пример на рис. 4).


\smallskip
{\it Доказательство необходимости в теореме.}
Пусть, напротив, некоторый $X$-граф с несамопересекающимися
циклами $K_1$ и $K_2$, имеющими ровно одну вершину $A$ перекрестья, $X$-планарен.
Цикл  $K_1$ делит плоскость на две части.
Лишь в вершине $A$ при движении по циклу $K_2$ мы переходим из одной части в другую,
а в остальных общих вершинах остаёмся в той же части. Противоречие.
QED

\smallskip
{\bf Задача 1.}
Если в $X$-графе найдутся два цикла с нечётным числом точек перекрестья, то найдутся и два цикла ровно с одной точкой перекрестья.


\smallskip
{\bf Задача 2} (для исследования).
Сформулируйте определения $\ast$-графа и $\ast$-планарности для графов, степень каждой вершины которых равна 4.
(Указание: полуребра, выходящие из каждой вершины, можно разбить на три пары или на две тройки...)
Найдите соответствующие критерии $\ast$-планарности.

Для одного из определений $\ast$-графа и $\ast$-планарности решение приведено в [F] и [FM].
Насколько мне известно, для других вариантов вопрос открыт!

\smallskip
В оставшейся части заметки приводится {\it доказательство достаточности в теореме.}

Цикл в графе называется {\it эйлеровым}, если он проходит по каждому {\it ребру}
ровно по одному разу.

Пусть в графе степени 4 задан эйлеров цикл.
Назовем две вершины $P$ и $Q$ {\it скрещивающимися}, если они идут вдоль эйлерова цикла
в порядке $PQPQ$ (а не $PPQQ$).
Назовем эйлеров цикл {\it удобным}, если вершины графа можно так разбить на два множества,
что никакие две вершины из одного множества не скрещиваются.

Идею доказательства можно пояснить на следующей красивой задаче.
\footnote{В дальнейшем она не используется формально, а служит лишь для
{\it мотивировки} понятий поворачивающего и сильно поворачивающего эйлерова цикла.
Поэтому задачу 3 и её решение можно пропустить.
}

\smallskip
{\bf Задача 3.} Связный граф степени 4 планарен тогда и только тогда, когда
в нём найдётся удобный эйлеров цикл.

\medskip
\centerline{\epsfbox{pict1.12}}
\centerline{Рисунок 5}

\smallskip
{\bf Решение.}
{\it Достаточность.}
Пусть вершины разбиты на два множества требуемым образом.
Построим несамопересекающийся цикл $S$ на плоскости, рёбра которого соответствуют всем
рёбрам удобного эйлерова цикла и идут в том же порядке.
Каждой вершине графа отвечают две вершины цикла $S$.
Соединим каждую пару вершин цикла $S$, отвечающую вершине {\it первого} множества
(из определения удобности),
несамопересекающейся ломаной {\it внутри} цикла $S$.
Эти ломаные можно взять попарно непересекающимися, поскольку соответствующие вершины
первого множества не скрещиваются.
Аналогично, соединим пары вершин цикла $S$, отвечающие вершинам {\it второго} множества,
не\-са\-мо\-пе\-ре\-се\-ка\-ю\-щи\-ми\-ся попарно непересекающимися ломаными {\it вне} цикла $S$.
Стянем в плоскости каждую проведённую ломаную в точку (рис. 5).
Получим вложение графа в плоскость.

{\it Необходимость.} Будем называть {\it гранями} связные куски, на которые распадается плоскость при разрезании по всем рёбрам плоского графа.
Начнём двигаться из некоторой вершины $A$, причём из каждой вершины будем выходить
по ребру, соседнему по грани с тем, по которому пришли.
По рёбрам разрешается проходить не более чем по одному разу.
Будем так двигаться, пока возможно.
Завершиться это движение может только в $A$.
Более того, мы вернёмся в $A$ по ребру, лежащему в одной грани с начальным
ребром (иначе мы ещё можем продолжить движение).
По окончании движения получим цикл $T'$, проходящий через каждую вершину
по парам рёбер, лежащих в одной грани (в том числе и через $A$, как показано ранее).
Если $T'$ не эйлеров, то из некоторой его вершины $B$ выходят рёбра, не принадлежащие $T'$.
Пройдём из $B$ по циклу $T'$, а затем продолжим движение так же, как при построении цикла $T'$.
В итоге получим больший цикл, также проходящий через каждую вершину по рёбрам, лежащим в одной грани,
и потому несамопересекающийся.
Таким образом можно увеличивать цикл, пока он не будет содержать все рёбра графа.

Докажем, что построенный эйлеров цикл удобен.
Расцепим этот цикл в каждой вершине, соединив в плоскости расцеплённые точки отрезками (рис. 6).
Объединение рёбер полученного графа, не совпадающих с проведёнными отрезками, разбивает плоскость на две части.
Включим в первое (второе) множество вершин исходного графа те вершины, которые отвечают проведённым отрезкам,
лежащим во внутренней (внешней) части. QED

\medskip
\centerline{\epsfbox{pict1.11}}
\centerline{Рисунок 6}

\smallskip
Пусть в $X$-графе имеется удобный эйлеров цикл $T$.
По нему можно построить вложение $X$-графа в плоскость, как в задаче 3.
При каком условии на $T$ полученное вложение будет <<$X$-вложением>>?
Ясно, что $T$ должен быть {\it поворачивающим}, т.е. в каждой вершине
последовательно проходить по полурёбрам из разных пар, выходящим из этой вершины.
Этого условия недостаточно.
Цикл $T$ должен быть {\it сильно поворачивающим}, т.е. при движении по циклу $T$ (в одном произвольном
направлении) мы должны входить в каждую вершину оба раза по полурёбрам одной и той же пары.

Достаточность в теореме вытекает из следующих трёх лемм.

\smallskip
{\bf Лемма о планарности.}
{\it Если в $X$-графе существует удобный сильно поворачивающий эйлеров цикл, то $X$-граф $X$-планарен.}

\smallskip
{\bf Лемма о поворачивающем цикле.} {\it В любом связном $X$-графе существует поворачивающий эйлеров цикл.}

\smallskip
{\bf Лемма о дополнительных свойствах.} {\it Пусть $X$-граф не содержит двух не\-са\-мо\-пе\-ре\-се\-ка\-ю\-щих\-ся
циклов без общих рёбер, имеющих ровно одну вершину перекрестья.
Тогда любой поворачивающий эйлеров цикл является удобным и сильно поворачивающим.}

\smallskip
{\it Доказательство леммы о планарности.}
(Это доказательство, кроме последнего абзаца, повторяет решение задачи 3.)
Построим несамопересекающийся цикл $S$ на плоскости, рёбра которого
соответствуют всем рёбрам удобного эйлерова цикла и идут в том же порядке.
Каждой вершине графа отвечают две вершины цикла $S$.

Соединим каждую пару вершин цикла $S$, отвечающую вершине {\it первого} множества (из определения удобности),
несамопересекающейся ломаной {\it внутри} цикла $S$.
Эти ломаные можно взять попарно непересекающимися, поскольку соответствующие
вершины первого множества не скрещиваются.
Аналогично, соединим пары вершин цикла $S$, отвечающие вершинам {\it второго}
множества, несамопересекающимися попарно непересекающимися ломаными {\it вне} цикла $S$.

Стянем в плоскости каждую проведённую ломаную в точку (рис. 5).
Получим вложение $X$-графа в плоскость.
Поскольку взятый эйлеров цикл сильно поворачивающий, при обходе вокруг каждой
вершины полурёбра из первой и из второй пары чередуются.
Значит, это вложение является <<$X$-вложением>>.
QED

\smallskip
{\it Доказательство леммы о поворачивающем цикле}
(аналогичное доказательству критерия существования эйлерова цикла).
Начнём двигаться из некоторой вершины $A$.
Будем двигаться, меняя пару полуребра в каждой вершине и проходя по рёбрам не более чем по одному разу, пока возможно.
Так как из каждой вершины выходит поровну полурёбер обоих пар, то завершиться это движение может только в $A$.
Более того, мы вернулись в $A$ по полуребру, не парному начальному полуребру
(иначе мы еще можем выйти по полуребру другой пары, а значит движение не завершилось).
По окончании движения получим поворачивающий цикл $F$.

Если $F$ не эйлеров, то из некоторой его вершины $B$ выходят рёбра, не принадлежащие $F$.
Пройдём из $B$ по циклу $F$, а затем продолжим движение так же, как при построении цикла $F$.
В итоге получим больший поворачивающий цикл.
Таким образом можно увеличивать цикл, пока он не будет содержать все рёбра $X$-графа.
QED

\smallskip
{\it Доказательство сильной поворачиваемости в лемме о дополнительных свойствах.}
Пусть, напротив, найдётся такая вершина $A$, что поворачивающий цикл $T$ входит в неё
первый раз по полуребру одной пары, а второй --- по полуребру другой пары.
Обозначим через $T_1$ цикл, соответствующий части цикла $T$ от первого
выхода из $A$ до второго входа в $A$.
Обозначим через $T_2$ цикл, соответствующий остатку цикла $T$.
Тогда $T_1$ и $T_2$ --- циклы, поворачивающие во всех вершинах, кроме вершины $A$,
являющейся для них вершиной перекрестья.
Если путь поворачивает в вершине, то эта вершина не может быть для него точкой перекрестья с другим путём.
Значит, остальные вершины не являются для циклов $T_1$ и $T_2$ точками перекрестья.
Но такой пары циклов нет по предположению.
QED

\smallskip
{\it Доказательство удобности в лемме о дополнительных свойствах.}
Построим вспомогательный граф $Q(T)$.
Вершины графа $Q(T)$ --- те же, что у исходного графа.
Две вершины графа $Q(T)$ соединяются ребром, если они скрещиваются
(относительно цикла $T$).

Удобность цикла $T$ равносильна двудольности графа $Q(T)$.
А это, в свою очередь, равносильно отсутствию в $Q(T)$ циклов нечётной длины.

Пусть поворачивающий цикл $T$ неудобен.
Возьмём наименьший цикл $A$ нечётной длины $2k+1\ge3$ в $Q(T)$.

Так как цикл $A$ минимален, то скрещиваются только соседние его вершины.
Пройдём по эйлерову циклу $T$, обозначая появляющиеся вершины цикла $A$ через $Y_1,Y_2,\ldots ,Y_{4k+2}$:
каждая из этих вершин получит два номера, отличающиеся на 3.
Более аккуратно, обозначим через

$\bullet$ $T$ отображение из $\{1,2,\ldots,n\}$ в множество вершин исходного графа,
выражающее порядок обхода эйлерова цикла $T$,

$\bullet$ $1\le y_1<y_2<\ldots<y_{4k+2}\le n$ числа, образами которых являются вершины цикла $A$,

$\bullet$ $Y_j$ вершину $T(y_j)$.

Без ограничения общности можно считать, что $Y_2=Y_5$, $Y_4=Y_7$, ..., $Y_{4k}=Y_1$, $Y_{4k+2}=Y_3$.

Положим $[i]:=(T(y_i),T(y_i+1),\ldots,T(y_{i+1}))$ для каждого $i\in\{1,2,\ldots,4k+1\}$.
(Неформально, это <<отрезок эйлерова цикла $T$, проходимый между $Y_i$ и $Y_{i+1}$>>.
Чтобы строго объяснить, что это такое, и нужно рассматривать $T$ как отображение.)

\smallskip
\centerline{\epsfbox{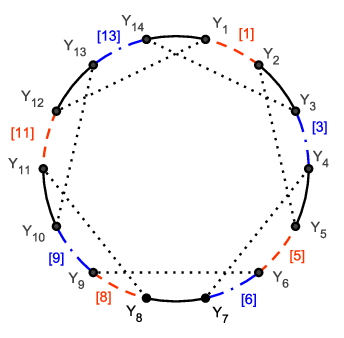}}
\centerline{Рисунок 7}
\smallskip

Возьмём первый цикл (рис.7):
$$[1],[5],\dots,[4l-7],[4l-3],\overline{[4l]},[4l+3],[4l+7],\dots,[4k-5],[4k-1],
\quad\mbox{где}\quad l:=[(k+1)/2].$$
Здесь путь $[4l]$ проходится против эйлерова цикла $T$, а остальные пути --- по циклу $T$.

(Вот неформальное объяснение этого построения для $k=3$, см. рис.7.
Пройдём по циклу $T$ из $Y_{12}=Y_1$ в $Y_2$, затем из $Y_2=Y_5$ в $Y_6$, затем против цикла $T$ из $Y_6=Y_9$ в $Y_8$, затем по циклу $T$ из $Y_8=Y_{11}$ в $Y_{12}$.)

Возьмём второй цикл, ``симметричный'' первому (рис.7):
$$[3][7],\dots,[4(k-l)-5],[4(k-l)-1],\overline{[4(k-l)+2]},[4(k-l)+5],[4(k-l)+9],\dots,[4k-3],[4k+1].$$
Здесь путь $[4(k-l)+2]$ проходится против цикла $T$, а остальные пути --- по циклу $T$.

(Вот неформальное объяснение этого построения для $k=3$, см. рис.7.
Пройдём по циклу $T$ из $Y_{14}=Y_3$ в $Y_4$, затем против цикла $T$ из $Y_7=Y_4$ в $Y_6$, затем по циклу $T$ из $Y_6=Y_9$ в $Y_{10}$, затем из $Y_{10}=Y_{13}$ в $Y_{14}$.)


У построенных двух циклов одна общая вершина $Y_{2k}=Y_{2k+3}$.
Первый цикл выходит из неё против цикла $T$, а входит по циклу $T$.
Как уже доказано, цикл $T$ сильно поворачивающий.
Поэтому первый цикл в обоих случаях идёт по полурёбрам одной пары.
Аналогично второй цикл проходит через эту вершину по полурёбрам другой пары.
Таким образом, у этих циклов единственная точка перекрестья.
Противоречие.
QED

\bigskip
\centerline{\sc Список литературы}
\medskip

[V] В.А. Васильев. {\it Инварианты и когомологии первого порядка для пространств
вложений самопересекающихся кривых в $\mathbb R^n$.} // Изв. РАН, сер. матем. 2005. 69:5.
С. 3-52.

[M1] В.О. Мантуров. {\it Доказательство гипотезы В. А. Васильева
о планарности сингулярных зацеплений.} // Изв. РАН, сер. матем. 2005. 69:5.
С. 169-178.

[M2] В.О. Мантуров. {\it Геометрия и комбинаторика виртуальных узлов.} Диссертация.  
\url{http://www.dissercat.com/content/geometriya-i-kombinatorika-virtualnykh-uzlov}

[M3] В. О. Мантуров. {\it Четырехвалентные графы с крестовой структурой.} // Математическое просвещение,  сер. 3, вып. 15. М.: МЦНМО. 2011. С. 128-137.

[F] T. Friesen. {\it A generalization of Vassiliev's planarity criterion.}

\url{http://arxiv.org/abs/1210.1539}

[FM] T. Friesen, V. Manturov. {\it Embeddings of $\ast$-graphs into 2-surfaces.}

\url{http://arxiv.org/abs/1212.5646}

\end{document}